\documentclass{article}

\usepackage{amssymb}
\usepackage{amsmath}
\usepackage[dvips]{graphicx}

\begin{document}


\noindent

\noindent

\noindent

\noindent

\begin{center} \textbf{ON AN  INTERIOR COMPACTNESS OF ONE
HOMOGENEOUS} \textbf{BOUNDARY -- VALUE PROBLEM}

\bigskip
\textbf{Rustamova Lamiya Aladdin}
\bigskip

 \emph{Institute of Applied Mathematics, Baku
State University, Z. Khalilov 23, AZ1148 Baku,
Azerbaijan,
 e-mail:\,rustamova@bsu.az}

\end{center}

\textbf{Abstract. } \emph{In the paper the conditions are obtained
providing existence and uniqueness  of the regular solution of the
boundary problem for class of the  second order homogeneous
operator-differential equation with singular coefficients. High term
of the equation contains the normal operator the spectrum of which
is contained in the certain sectors. }

 \emph{ Further, it is proved
the theorem of internal compactness of space of regular solutions of
the considered  problem.}


\bigskip

 We consider the boundary problem for homogeneous operator-differential equation in separable Hilbert space$H$

\begin{equation} \label{GrindEQ__1_}
-\frac{d^{2} u}{dt^{2} } +\rho (t)A^{2} u+A_{0} \frac{d^{2} u}{dt^{2} } +A_{1} \frac{du}{dt} +A_{2} u=0,\, \, \, \, \, \,
\end{equation}

\begin{equation} \label{GrindEQ__2_}
\, \, \, u(0)=\varphi \, \, \, \, ,\, \, \, \, \, \,
\end{equation}
where  $\varphi \in H_{{\raise0.5ex\hbox{$\scriptstyle 3$}
\kern-0.1em/\kern-0.15em \lower0.25ex\hbox{$\scriptstyle 2$}}}
 $, $u(t)\in W_{2}^{2} (R_{+} ;H)$, $\rho (t)$  is the form of

\noindent                                       $\rho (t)=\left\{\begin{array}{l} {\alpha ^{2} ,\, \, \, \, t\in (0;1),} \\ {\beta ^{2} ,\, \, \, \, t\in (1;\infty )\, \, ,} \end{array}\right. $                                                                                                                                                   moreover  $\alpha >0,\, \, \, \beta >0$,  operator coefficients  $A$ and  $A_{j} \, \, \, \, (j=0,1,2)$ satisfy the following conditions

\noindent 1) $A$ is normal, with quite continuous inverse $A^{-1} $ operator, spectrum of which is contained in a corner sector,

\[S_{\varepsilon } =\left\{\lambda :\left|\arg \lambda \right|\le \varepsilon \, \right\}\, \, \, \, \, ,\, \, \, \, \, 0\le \varepsilon <{\raise0.7ex\hbox{$ \pi  $}\!\mathord{\left/{\vphantom{\pi  2}}\right.\kern-\nulldelimiterspace}\!\lower0.7ex\hbox{$ 2 $}} ;\]
2) Operators  $B_{j} =A_{j} A^{-j} \, \, \, \, (j=0,1,2,)$are bounded in  $H$.

Denote by   $L_{2} (R_{+} ;H)$ a Hilbert space of the vector-functions $f(t)$ with values from $H$, measurable and integrable by square-law with norm

\[\left\| f\right\| _{L_{2} } =\left(\int _{0}^{+\infty }\left\| f(t)\right\| ^{2} dt \right)
^{{\raise0.5ex\hbox{$\scriptstyle 1$} \kern-0.1em/\kern-0.15em
\lower0.25ex\hbox{$\scriptstyle 2$}}}
 .\]
Further we introduce the  space e.g. [1]

\[W_{2}^{2} (R_{+} ;H)=\left\{u(t):\, \, \, \, u'',A^{2} u\in L_{2} (R_{+} ;H)\, \, \, ,\, \left\| u\right\| _{W_{2}^{2} }
=\left(\left\| u''\right\| _{_{L_{2}^{} } }^{2} +\left\| A^{2}
u\right\| _{_{L_{2,}^{} } }^{2} \right )
^{{\raise0.5ex\hbox{$\scriptstyle 1$} \kern-0.1em/\kern-0.15em
\lower0.25ex\hbox{$\scriptstyle 2$}}}  \, \right\}.\] Then from the
trace theorem it  results that

\[\mathop{W_{2}^{2} }\limits^{o} (R_{+} ;H;\left\{0\right\})=\left\{u\, \, \, \left|\, u\in W_{2}^{2} (R_{+} ;H)\, ,\, u(0)=0\right. \right\}.\]
    \textbf{Definition. }If for any\textbf{ }  $\varphi \in H_{{\raise0.5ex\hbox{$\scriptstyle 3$}
\kern-0.1em/\kern-0.15em \lower0.25ex\hbox{$\scriptstyle 2$}}}$
there exists the vector-function  $u(t)$ which satisfies
\eqref{GrindEQ__1_}, and  boundary condition \eqref{GrindEQ__2_} in
the sense

\[\mathop{\lim }\limits_{t\to +0} \, \left\| u(t)-\varphi \right\| _{{\raise0.5ex\hbox{$\scriptstyle 3$}
\kern-0.1em/\kern-0.15em \lower0.25ex\hbox{$\scriptstyle 2$}}} =0\]
also the  estimation takes place

\[\left\| u\right\| _{W_{2}^{2} (R_{+} ;H)}^{} \le const\left\| \varphi \right\| _
{_{{\raise0.5ex\hbox{$\scriptstyle 3$} \kern-0.1em/\kern-0.15em
\lower0.25ex\hbox{$\scriptstyle 2$}}} }^{}   ,\] then $u(t)$ is
called a regular solution of the problem
\eqref{GrindEQ__1_},\eqref{GrindEQ__2_}, and the problem
\eqref{GrindEQ__1_},\eqref{GrindEQ__2_}  is called regular solvable.

 We shall note that when the equations are  not  homogeneous and $A_{0} =1,\, \, \, \, \alpha =\beta $=1 this problem is investigated in [2], when  $A=A^{*} \ge cE,\, \, c>0$, at $A_{0} =1,\, \, \, \, \alpha \ne \beta $ in  [3]. When the equation is non homogeneous boundary problem \eqref{GrindEQ__1_},\eqref{GrindEQ__2_} is investigated [4] and resolvability of the equation \eqref{GrindEQ__1_} on all axis it is considered in  [5] .

  First we shall consider the  problem

\begin{equation} \label{GrindEQ__3_}
P_{o} u=-\frac{d^{2} u}{dt^{2} } +\rho (t)A^{2} u=0\, \, \, \, ,\, \, \, \,
\end{equation}

\begin{equation} \label{GrindEQ__4_}
u(0)\, =\varphi \, \, \, \, \, \, \, \, \, .\, \, \, \, \, \, \, \,
\end{equation}
Let's seek  the regular solution of the  problem \eqref{GrindEQ__3_}, \eqref{GrindEQ__4_} in the form

\[u_{0} (t)=\left\{\begin{array}{l} {e^{-\alpha tA} c_{1} +e^{-\alpha (1-t)A} c_{2} \, \, \, ,\, \, \, \, \, \, \, \, \, \, \, \, \, \, \, \, \, t\in (0;\, \, 1),} \\ {e^{\beta A(1-t)} c_{3} \, \, \, ,\, \, \, \, \, \, \, \, \, \, \, \, \, \, \, \, \, \, \, \, \, \, \, \, \, \, \, \, \, \, \, \, \, \, \, \, \, \, \, \, t\in (1;\infty ),} \end{array}\right. \]
where  $c_{1} \, ,c_{2} \, ,c_{3} $ - are unknown elements from $H
_{{\raise0.5ex\hbox{$\scriptstyle 3$} \kern-0.1em/\kern-0.15em
\lower0.25ex\hbox{$\scriptstyle 2$}}}  $. From  the  condition
\eqref{GrindEQ__4_} and inclusion  $u_{0} (t)\in W_{2}^{2} (R_{+}
;H)$ it is obtained  the following system of the equations
relatively  $c_{1} \, ,c_{2} $  and  $\, c_{3} $ :

\[\left\{\begin{array}{l} {c_{1} +e^{-\alpha A} c_{2} =\varphi ,} \\ {e^{-\alpha A} c_{1} +c_{2} \, \, =c_{3} ,\, \, \, \, } \\ {\, -\alpha Ae^{-\alpha A} c_{1} +\alpha Ac_{2} \, \, \, \, =-\beta Ac_{3} \, \, \, ,\, \, \, \, \, \, \, \, \, \, \, \, \, \, \, \, \, } \end{array}\right. \]
or

\[\left\{\begin{array}{l} {c_{1} +e^{-\alpha A} c_{2} +0\cdot c_{3} =\varphi \, \, ,} \\ {e^{-\alpha A} c_{1} +c_{2} \, \, -c_{3} =0\, \, ,} \\ {\, -\alpha e^{-\alpha A} c_{1} +\alpha c_{2} \, \, \, +\beta c_{3} \, \, =0\, \, ,\, \, \, \, \, \, \, \, \, \, \, \, \, \, \, \, \, \, \, \, \, \, \, \, } \end{array}\right. \]
or in an operational

\[\Delta _{0} (A)c=\tilde{\varphi },\]
where

\[\Delta _{0} (A)=\left|\begin{array}{l} {E\, \, \, \, \, \, \, \, \, \, \, \, \, \, \, \, e^{-\alpha A} \, \, \, \, \, \, \, \, \, \, \, \, 0} \\ {-e^{-\alpha A} \, \, \, \, \, \, \, E\, \, \, \, \, \, \, \, \, \, \, \, \, \, \, -E} \\ {-\alpha \, e^{-\alpha A} \, \, \, \, \, \, \, \, \alpha E\, \, \, \, \, \, \, \, \, \, \, \, \, \beta E} \end{array}\right|,\]

\[c=\left(\begin{array}{l} {c_{1} } \\ {c_{2} } \\ {c_{3} } \end{array}\right)\, ,       \tilde{\varphi }=\left(\begin{array}{l} {\varphi } \\ {0} \\ {0} \end{array}\right).\]
As we have shown $\Delta _{0}^{} (A)$ that it is  $H^{3} =H\times
H\times H$(see . [4]), therefore, we shall unequivocally define
$c_{1} \, ,c_{2} $ and $\, c_{3} $. They belong
$H_{{\raise0.5ex\hbox{$\scriptstyle 3$} \kern-0.1em/\kern-0.15em
\lower0.25ex\hbox{$\scriptstyle 2$}}}  $, as $\varphi \in
H_{{\raise0.5ex\hbox{$\scriptstyle 3$} \kern-0.1em/\kern-0.15em
\lower0.25ex\hbox{$\scriptstyle 2$}}}  $. It is obvious, that

\[\left\| u_{0} (t)\right\| _{W_{2}^{2} (R_{+} ;H)}^{} \le const\left\| \varphi \right\|
 _{{\raise0.5ex\hbox{$\scriptstyle 3$} \kern-0.1em/\kern-0.15em
\lower0.25ex\hbox{$\scriptstyle 2$}}}   .\] Now we  consider the
boundary problem \eqref{GrindEQ__1_},\eqref{GrindEQ__2_}. For this
purpose we take

\[u(t)=\vartheta (t)+u_{0} (t).\]
Then we get the following equation

\[\begin{array}{l} {-\frac{d^{2} \vartheta (t)}{dt^{2} } +\rho (t)A^{2} \vartheta (t)+A_{0}
 \frac{d^{2} \vartheta (t)}{dt^{2} } +A_{1} \frac{d\vartheta (t)}{dt} +A_{2} \vartheta (t)-} \\
 \\ {-\, \,
  \frac{d^{2} u_{0} (t)}{dt^{2} } +\rho (t)A^{2} u_{0} (t)+A_{0} \frac{d^{2} u_{0} (t)}{dt^{2} }
   +A_{1} \frac{du_{0} (t)}{dt} +A_{2} u_{0} (t)\, \, =0,\, \, \, \, \, \, \, \, \, \, \, \, \, \, \, \, \, \, \, \, \, \, \, \, \, \, } \\ {\, \, \, \, \, \, \, \, \, \, \, \, \, \, \, \, \, \, \, \, \, \, \, \, \, \, \, \, \, \, \, \, \, \, \, \, \, \, \, \, \, \, \, \, \, \, \, \, \, \, \, \, \, \, \, \, \, \, \, \, \, \, \, \, \, \, \, \, \, \, } \end{array}\]

\[\vartheta (0)=0.\]
As      $$-\frac{d^{2} u_{0} (t)}{dt^{2} } +\rho (t)A^{2} u_{0}
(t)=0\, \, \, \, ,\, \, \, \, $$

\[\begin{array}{l} {-\frac{d^{2} \vartheta (t)}{dt^{2} } +\rho (t)A^{2} \vartheta (t)+A_{0} \frac{d^{2} \vartheta (t)}{dt^{2} } +A_{1} \frac{d\vartheta (t)}{dt} +A_{2} \vartheta (t)=g(t)\, \, \, \, \, \, \, \, \, \, \, \, \, \, \, \, \, ,\, \, \, \, \, \, \, \, \, \, \, \, \, \, \, \, \, \, \, \, } \\ {\, \, \, \, \, \, \, \, \, \, \, \, \, \, \, \, \, \, \, \, \, \, \, \, \, \, \, \, \, \, \, \, \, \, \, \, \, \, \, \, \, \, \, \, \, \, \, \, \, \, \, \, \, \, \, \, \, \, \, \, \, \, \, \, \, \, \, \, \, \, } \end{array}\]

\[\vartheta (0)=0\, \, \, \, \, \, \, \, \, ,\]
where

\[\begin{array}{l} {g(t)=\, -A_{0} \frac{d^{2} u_{0} (t)}{dt^{2} } -A_{1} \frac{du_{0} (t)}{dt} -A_{2} u_{0} (t)\, \, \, .} \\ {\, \, \, \, \, \, \, \, \, \, \, \, \, \, \, \, \, \, \, \, \, \, \, \, \, \, \, \, \, \, \, \, \, \, \, \, \, \, \, \, \, \, \, \, \, \, \, \, \, \, \, \, \, \, \, \, \, \, \, \, \, \, \, \, \, \, \, \, \, \, } \end{array}\]
As operators  $A_{j} A^{-j} \, \, (j=\overline{0,2})$ are bounded

\[
\begin{array}{l}
 \left\| {g(t)} \right\|_{L_2 (R_ +  ;H)}  \le \left\| {A_0 } \right\|\left\| {\frac{{d^2 u_0 }}{{dt^2 }}} \right\|_{L_2 (R_ +  ;H)}  +  \\
  \\
  + \left\| {A_1 A^{ - 1} } \right\|\left\| {A\frac{{du_0 }}{{dt}}} \right\|_{L_2 (R_ +  ;H)}  + \left\| {A_2 A^{ - 2} } \right\|\left\| {A^2 u_0 } \right\|_{L_2 (R_ +  ;H)} . \\
  \\
 \end{array}
\]

Applying , the theorem of intermediate derivatives, [1] we have

\[\left\| g(t)\right\| _{L_{2} (R_{+} ;H)} \le const\left\| u_{0} \, \right\| _{W_{2}^{2} (R_{+} ;H)} \,
 \le const\left\| \varphi \right\| _{{\raise0.5ex\hbox{$\scriptstyle 3$} \kern-0.1em/\kern-0.15em
\lower0.25ex\hbox{$\scriptstyle 2$}}}  \, \, \, .\] Thus, we have
reduced a boundary problem \eqref{GrindEQ__1_},\eqref{GrindEQ__2_}
 to the non homogeneous boundary problem with a zero boundary condition. Thus  following theorem is valid.

 \textbf{Theorem 1. }Let the operator $A$ satisfies the condition  1), but operators
  $B_{j} =A_{j} A^{-j} \, \, \, \, (j=0,1)$ satisfy condition    2) , moreover

\noindent

\noindent where  numbers $c_{j} (\varepsilon ;\alpha ;\beta )$ are defined as

\[
c_0 (\varepsilon ;\alpha ;\beta ) = \frac{1}{{\min (\alpha ^2 ;\beta
^2 )}}\left\{ \begin{array}{l}
 1,\,\,\,\,\,\,\,\,\,\,\,\,\,\,\,\,\,\,0 \le \varepsilon  < {\raise0.7ex\hbox{$\pi $} \!\mathord{\left/
 {\vphantom {\pi  {4}}}\right.\kern-\nulldelimiterspace}
\!\lower0.7ex\hbox{${4}$}}, \\
  \\
 \frac{1}{{\sqrt 2 \cos \varepsilon }},\,\,\,\,\,\, {\raise0.7ex\hbox{$\pi $}\!\mathord {\left/
 {\vphantom {\pi  4}}\right.\kern-\nulldelimiterspace}
\!\lower0.7ex\hbox{$4$}} \le \varepsilon  < {\raise0.7ex\hbox{$\pi
$} \!\mathord{\left/
 {\vphantom {\pi  2}}\right.\kern-\nulldelimiterspace}
\!\lower0.7ex\hbox{$2$}}. \\
 \end{array} \right.
\]

\[c_{1} (\varepsilon ;\alpha ;\beta )=\frac{1}{2\, \cos \varepsilon \, \, \min (\alpha ;\beta )} \, ,\, \, \, \, \, \, \, \, \, \, \, \, \, \, 0\, \, \, \, \le \varepsilon <\, {\raise0.7ex\hbox{$ \pi  $}\!\mathord{\left/{\vphantom{\pi  2}}\right.\kern-\nulldelimiterspace}\!\lower0.7ex\hbox{$ 2 $}} ,\]

\[
c_2 (\varepsilon ;\alpha ;\beta ) = \frac{{\max (\alpha ;\beta
)}}{{\min (\alpha ^2 ;\beta ^2 )}}\left\{ \begin{array}{l}
 1,\,\,\,\,\,\,\,\,\,\,\,\,\,\,\,\,\,\,0 \le \varepsilon  \le {\raise0.7ex\hbox{$\pi $} \!\mathord{\left/
 {\vphantom {\pi  {4}}}\right.\kern-\nulldelimiterspace}
\!\lower0.7ex\hbox{${4}$}}, \\
  \\
 \frac{1}{{\sqrt 2 \cos \varepsilon }},\,\,\,\,\,\,\,{\raise0.7ex\hbox{$\pi $} \!\mathord{\left/
 {\vphantom {\pi  4}}\right.\kern-\nulldelimiterspace}
\!\lower0.7ex\hbox{$4$}} \le \varepsilon  < {\raise0.7ex\hbox{$\pi
$} \!\mathord{\left/
 {\vphantom {\pi  2}}\right.\kern-\nulldelimiterspace}
\!\lower0.7ex\hbox{$2$}}. \\
 \end{array} \right.
\]
Then boundary problem \eqref{GrindEQ__1_}, \eqref{GrindEQ__2_} is regularly solvable.

 Now we shall study one property of homogeneous regular solutions. Let numbers $a,a_{1} ,b_{1} ,b$  be such, that

\[0<a<a_{1} <b_{1} <b<\infty .\]
Denote by  $N(P)$ the  space of regular solutions  of the  boundary  problem \eqref{GrindEQ__1_}, \eqref{GrindEQ__2_}. It is obvious, that $N(P)$- linear full subspace  in $W_{2}^{2} (R_{+} ;H)$. Really, if $u_{n} (t)\to u(t)$ in $W_{2}^{2} (R_{+} ;H)$ and$P\left({d \mathord{\left/{\vphantom{d dt}}\right.\kern-\nulldelimiterspace} dt} \right)u_{n} (t)=0\, \, \, \, \, (u_{n} (t)\in N(P))$,

\[\left\| P\left({d \mathord{\left/{\vphantom{d dt}}\right.\kern-\nulldelimiterspace} dt} \right)\left(u(t)-u_{n} (t)\right)\right\| \le const\left\| u(t)-u_{n} (t)\right\| \, \, \to 0.\]
Then $\left\| P\left({d \mathord{\left/{\vphantom{d dt}}\right.\kern-\nulldelimiterspace} dt} \right)u(t)\right\| =0,$ i.e. $P\left({d \mathord{\left/{\vphantom{d dt}}\right.\kern-\nulldelimiterspace} dt} \right)u(t)=0,$hence $u(t)\in N(P)$.

\noindent It is obvious, that $N(P)\subset W_{2}^{1} (R_{+} ;H)$.

 \textbf{ Definition 2.} If $a,a_{1} ,b_{1} ,b$ are such, as $0<a<a_{1} <b_{1} <b<\infty $$M>0$ , the set $\left\{u|u\in N(P),\, \, \left\| u\right\| _{W_{2}^{1} \left((a,b);H\right)} \le M\right\}$ is compact on norm $\left\| u\right\| _{W_{2}^{1} ((a_{1} ,b_{1} );H)} $we say  speak, that space of regular solutions  of  the  problem \eqref{GrindEQ__1_}, \eqref{GrindEQ__2_} is internally compact.

  We note, that definition of internal compactness for the first time has entered P.D.Laks [6]. At different situations of interior compactness of solutions the considered works [7, 8]. Following P.D.Laks's [6] work, we have entered concept of interior  compactness of the solutions of the homogeneous equations.

  Takes place

 \textbf{Theorem 2.} A condition of the theorem 1 let satisfied. Then the space of regular solutions  of a problem \eqref{GrindEQ__1_}, \eqref{GrindEQ__2_} is internally compact.

 \textbf{Proof.} Let $0<a<a_{1} <b_{1} <b<\infty $ and scalar function $\varphi (t)\in C_{0}^{\infty } (a,\, b)$, is such, that

\[\varphi (t)=\left\{\begin{array}{l} {1,\, \, \, \, \, \, \, \, \, \, \, \, \, \, \, \, \, t\in (a_{1} ,b_{1} ),} \\ {0,\, \, \, \, \, \, \, \, \, \, \, \, \, \, \, \, \, t\ge b,\, \, \, t\le a.} \end{array}\right. \]
Then it is obvious, that for vector functions $\varphi (t)\, u\, (t)\, \in \mathop{W_{2}^{2} }\limits^{o} (R_{+} ;H;0)$ and $u(t)=0$ at $t\le a$ and $t\ge b$.

  As we have proved, that at performance of conditions of the theorem the following inequality takes place:

\[\left\| P\left({d \mathord{\left/{\vphantom{d dt}}\right.\kern-\nulldelimiterspace} dt} \right)\varphi u\right\| _{L_{2} (R_{+} ;H)} \ge const\left\| \varphi u\right\| _{W_{2}^{2} (R_{+} ;H)} \, \, \]
For all $u\, \, \in \mathop{W_{2}^{2} }\limits^{o} (R_{+} ;H;0)$ (see [4], the theorem 2).

\noindent From here we have:

\[\begin{array}{l} {\left\| -\frac{d^{2} \varphi (t)u(t)}{dt^{2} } +\rho (t)\varphi (t)A^{2} u(t)+A_{0} \frac{d^{2} \varphi (t)u(t)}{dt^{2} } +A_{1} \frac{d\varphi (t)u(t)}{dt} +A_{2} \varphi (t)u(t)\right\| _{L_{2} (R_{+} ;H)} \ge } \\ {\ge const\left\| \varphi u\right\| _{W_{2}^{2} (R_{+} ;H)} \, \, \, \, .\, } \end{array}\]
As

\[
\begin{array}{l}
 P\left( {{d \mathord{\left/
 {\vphantom {d {dt}}} \right.
 \kern-\nulldelimiterspace} {dt}}} \right)\varphi (t)u(t) =  - \frac{{d^2 \varphi (t)u(t)}}{{dt^2 }} + \rho (t)\varphi (t)A^2 u(t) + A_0 \frac{{d^2 \varphi (t)u(t)}}{{dt^2 }} +  \\
  \\
  + A_1 \frac{{d\varphi (t)u(t)}}{{dt}} + A_2 \varphi (t)u(t) =  - \frac{{d^2 \varphi (t)}}{{dt^2 }}u(t) - 2\frac{{d\varphi (t)}}{{dt}}\frac{{du(t)}}{{dt}} - \varphi (t)\frac{{d^2 u(t)}}{{dt^2 }} +  \\
  \\
  + \rho (t)\varphi (t)A^2 u(t) + A_0 \left( {\frac{{d^2 \varphi (t)}}{{dt^2 }}u(t) + 2\frac{{d\varphi (t)}}{{dt}}\frac{{du(t)}}{{dt}} + \frac{{d^2 u(t)}}{{dt^2 }}\varphi (t)} \right) +  \\
  \\
  + A_1 \left( {\frac{{d\varphi (t)}}{{dt}}u(t) + \frac{{du(t)}}{{dt}}\varphi (t)} \right) + A_2 \varphi (t)u(t) =  \\
 \end{array}
\]

\[
\begin{array}{l}
  = \varphi (t)\left( { - \frac{{d^2 u(t)}}{{dt^2 }} + \rho (t)A^2 u(t) + A_0 \frac{{d^2 u(t)}}{{dt^2 }} + A_1 \frac{{du(t)}}{{dt}} + A_2 u(t)} \right) +  \\
  \\
  + \left( { - \frac{{d^2 \varphi (t)}}{{dt^2 }}u(t) - 2\frac{{d\varphi (t)}}{{dt}}\frac{{du(t)}}{{dt}} + A_0 \left( {\frac{{d^2 \varphi (t)}}{{dt^2 }}u(t) + 2\frac{{d\varphi (t)}}{{dt}}\frac{{du(t)}}{{dt}}} \right) + A_1 \frac{{d\varphi (t)}}{{dt}}u(t)} \right). \\
 \end{array}
\]
As $\varphi (t)=0$ at $t\ge b$, $t\le a$ and $u(t)$ - the regular decision,

\[P\left({d \mathord{\left/{\vphantom{d dt}}\right.\kern-\nulldelimiterspace} dt} \right)u(t)=-\frac{d^{2} u(t)}{dt^{2} } +\rho (t)A^{2} u(t)+A_{0} \frac{d^{2} u(t)}{dt^{2} } +A_{1} \frac{du(t)}{dt} +A_{2} u(t)=0\, \, \, \]
And

\[
\begin{array}{l}
 \left\| {P\left( {{d \mathord{\left/
 {\vphantom {d {dt}}} \right.
 \kern-\nulldelimiterspace} {dt}}} \right)u(t)} \right\|_{L_2 (R_ +  ;H)}  = \left\| { - \frac{{d^2 \varphi (t)}}{{dt^2 }}u(t) - 2\frac{{d\varphi (t)}}{{dt}}\frac{{du(t)}}{{dt}} + } \right. \\
  \\
  + \left. {A_0 \left( {\frac{{d^2 \varphi (t)}}{{dt^2 }}u(t) + 2\frac{{d\varphi (t)}}{{dt}}\frac{{du(t)}}{{dt}}} \right) + A_1 \frac{{d\varphi (t)}}{{dt}}u(t)} \right\|_{L_2 (R_ +  ;H)}  \ge  \\
  \\
  \ge const\left\| {\varphi u} \right\|_{W_2^2 (R_ +  ;H)}  = const\left\| {\varphi u} \right\|_{W_2^2 ((a;b);H)}  \ge  \\
  \\
  \ge const\left\| {\varphi u} \right\|_{W_2^2 ((a_1 ;b_1 );H)}  = const\left\| u \right\|_{W_2^2 ((a_1 ;b_1 );H)} . \\
 \end{array}
\]
Thus,

\[
\begin{array}{l}
 \left\| u \right\|_{W_2^2 ((a_1 ;b_1 );H)}  \le \left\| { - \frac{{d^2 \varphi (t)}}{{dt^2 }}u(t) - 2\frac{{d\varphi (t)}}{{dt}}\frac{{du(t)}}{{dt}}} \right. +  \\
  \\
  + \left. {A_0 \left( {\frac{{d^2 \varphi (t)}}{{dt^2 }}u(t) + 2\frac{{d\varphi (t)}}{{dt}}\frac{{du(t)}}{{dt}}} \right) + A_1 \frac{{d\varphi (t)}}{{dt}}u(t)} \right\|_{L_2 (R_ +  ;H)}  \le  \\
  \\
  \le const\left\| {u(t)} \right\|_{L_2 ((a;b);H)}  + const\left\| {\frac{{du(t)}}{{dt}}} \right\|_{L_2 ((a;b);H)}  + const\left\| {A_0 u(t)} \right\|_{L_2 ((a;b);H)}  +  \\
  \\
  + const\left\| {A_0 \frac{{du(t)}}{{dt}}} \right\|_{L_2 ((a;b);H)}  + const\left\| {A_1 u(t)} \right\|_{L_2 ((a;b);H)}  \le  \\
  \\
  \le const\left( {\left\| u \right\|_{L_2^{} ((a;b);H)}  + \left\| {\frac{{du(t)}}{{dt}}} \right\|_{L_2^{} ((a;b);H)}  + \left\| {A_0 } \right\|\left\| u \right\|_{L_2^{} ((a;b);H)}  + } \right. \\
  \\
  + \left. {\left\| {A_0 } \right\|\left\| {\frac{{du(t)}}{{dt}}} \right\|_{L_2^{\scriptstyle  \hfill \atop
  \scriptstyle 2 \hfill} ((a;b);H)}  + \left\| {A_1 A^{ - 1} } \right\|\left\| {Au} \right\|_{L_2^{} ((a;b);H)} } \right) \le const\left\| u \right\|_{W_2^1 ((a;b);H)}  \le M. \\
 \end{array}
\]

Hence, for anyone $0<a<a_{1} <b_{1} <b$we have:

\[\left\| u\right\| _{W_{2}^{2} ((a_{1} ;b_{1} );H)} \le const\left\| u\right\| _{W_{2}^{1} ((a;b);H)} \, \, \le M.\]
As $u\in N(P)$, the set $\left\{u|\left\| u\right\| _{W_{2}^{2} \left((a_{1} ;b_{1} );H\right)} ,u\in N(P)\right\}$ is limited. From quite continuity of the operator $A^{-1} $ follows, that the space $W_{2}^{2} \left(\left(a_{1} ,b_{1} \right);H\right)$is enclosed in space $W_{2}^{1} \left(\left(a_{1} ,b_{1} \right);H\right)$ compactly, i.e.

\[W_{2}^{2} \left(\left(a_{1} ,b_{1} \right);H\right)\subset W_{2}^{1} \left(\left(a_{1} ,b_{1} \right);H\right)\]
compactly. Hence, $N(P)$- compact set of century $W_{2}^{1} \left(\left(a_{1} ,b_{1} \right);H\right)$. Thus, we have proved internal compactness of decisions of a problem \eqref{GrindEQ__1_}, \eqref{GrindEQ__2_}. The theorem is proved.

\newpage
\begin{center}\textbf{References}\end{center}

\begin{enumerate}

\item  Lions ZH.-L., Madzhenes E. Non homogeneous boundary -- value  problems and their . M., Mir , 1971,371 p.

\item  Mirzoev S.S. The problems to theory to solvability of the boundaryvalue  problems for operator -differential of the equations in Hilbert space and spectral problems connected with them. The Thesis on competition dissert. of uch.step.dokt.fiz.-mat. Sciences.Baku 1993, 229 p.

\item  Mirzoev S.S., Aliev A.R. On  one boundary-value   problem for operator-differential equations of the second order with discontinuous coefficient. The Works of the Institute Mat. and Mech. AN ., VI (XIV),Baku, 1997, p.p.117-121.

\item  Mirzoyev S.S., Rustamova L.A. On solvability of on boundary-value problem for operator - differential equations of the second order with discontinuous coefficient //An International  Journal of Applied and Computational Mathematics, 2006, v. 5, ¹2, p.191-200.

\item  Rustamova L.A. On  regular solvability of one class operator - differential equations of the second order.  BGU , Vestnik BGU ,  2005,¹1, p.p.43-51.

\item  Lax P.D. A Phragmen-Lindelþf theorem in harmonic analysis and its application to some
          questions in the theory of elliptic equations // Comm.Pure Appl.Math., 10 (1957), pp.361-
          389.

\item  P. Koosis, \textit{Interior compact spaces of functions on a half-line, }Comm. Pure Appl. Math. 10 (1957), pp. 583-615.

\item  D. Baranov , Interior-compact subspaces and differentiation in model subspaces,  Journal   of   Mathematical Sciences, Volume 139, Number 2 , 2006, pp. 6369-6373

\end{enumerate}

\end{document}